\documentclass[a4paper, 12pt]{article}
\usepackage{amsmath, amsthm}
\usepackage[bookmarksnumbered]{hyperref}
\usepackage{amsrefs}

\newtheorem{theorem}{Theorem}[section]
\newtheorem{lemma}[theorem]{Lemma}

\numberwithin{equation}{section}

\author{G\'abor Braun\\
Alfr\'ed R\'enyi Institute of Mathematics\\
Hungarian Academy of Sciences\\
Budapest\\
Re\'altanoda u 13--15\\
H-1053}
\title{A proof of Higgins's conjecture\thanks{MSC 2000 Primary: 20E06,
    Secondary: 20L05}}
\date{}

\begin{document}
\maketitle

\begin{abstract}
  Let $\Theta\colon G=\prod^*_{\lambda \in \Lambda} G_\lambda \to B=\prod^*_{\lambda \in \Lambda} B_\lambda$ be a group homomorphism
  between free products of groups such that $G_\lambda\Theta=B_\lambda$ for all $\lambda \in \Lambda$. Let $H
  \subseteq G$ be a subgroup such that $H\Theta=B$.  Then $H = \prod^*_{\lambda \in \Lambda} H_\lambda$ such that
  $H_\lambda\Theta=B_\lambda$ and $H_\lambda = \prod^* (H \cap G_\lambda^{\beta_{\lambda,\mu}}) \ast F_\lambda$ where $F_\lambda$ is free.
\end{abstract}

\section{Introduction}
\label{sec:introduction}

Recall that the free product of groups $G_\lambda$ is the group $\prod^*_{\lambda \in \Lambda} G_\lambda$
generated by the $G_\lambda$ in which every relation follows from
group identities. In other words, free product is the same as the coproduct in the
category of groups. We use $\prod^*$ or $\ast$ to denote free products.

Let $H^x:=x^{-1}Hx$ denote a conjugate of a group $H$.
There are two main theorems about subgroups of free products of groups:
\begin{theorem}[Kuro\v s's Theorem]
  \label{th:Kurosh}
  Let $H \subseteq \prod^*_{\lambda \in \Lambda} G_\lambda$ be a subgroup of a free product. Then $H$ has a
  free decomposition $H=\prod^*_{\lambda \in \Lambda, x_\lambda} (H \cap G_\lambda^{x_\lambda}) \ast F$ where for each
  $\lambda$ the $x_\lambda$ runs
  through a suitable set of representatives of double cosets $G_\lambda x H$ such
  that $G_\lambda H$ is represented by $1$.  Moreover, $F$ is free.
\end{theorem}

\begin{theorem}[Higgins's Theorem]
  \label{th:Grushko-Higgins}
  Let $\Theta\colon G=\prod^*_{\lambda \in \Lambda} G_\lambda \to B=\prod^*_{\lambda \in \Lambda} B_\lambda$ be a group homomorphism such
  that $G_\lambda\Theta=B_\lambda$ for all $\lambda \in \Lambda$. Let $H \subseteq G$ be a subgroup such that $H\Theta=B$.
  Then there are groups $H_\lambda$ such that $H = \prod^*_{\lambda \in \Lambda} H_\lambda$ and $H_\lambda\Theta=B_\lambda$.
\end{theorem}

Higgins proved the above theorems in \cite[Chapter~14]{MR48:6288} using
groupoids.  These proofs are similar and Higgins conjectured that they can be
united to a single proof of a common generalization of the two theorems.

However, P.~R.~Heath and P.~Nickolas showed in~\cite{MR97f:20032} that there
are difficulties in generalizing the proof and, in particular, E.~T.~Ordman's
proof of Higgins's conjecture in~\cite{MR43:2102} is incorrect.

Nevertheless, we give a simple proof of Higgins's conjecture in this paper
using the two theorems above.

\begin{theorem}[Higgins's conjecture]
  \label{th:Higgins-conjecture}
  Let $\Theta\colon G=\prod^*_{\lambda \in \Lambda} G_\lambda \to B=\prod^*_{\lambda \in \Lambda} B_\lambda$ be a group homomorphism such
  that $G_\lambda\Theta=B_\lambda$ for all $\lambda \in \Lambda$. Let $H \subseteq G$ be a subgroup such that $H\Theta=B$.
  Then $H = \prod^*_{\lambda \in \Lambda} H_\lambda$ such that $H_\lambda\Theta=B_\lambda$ where $H_\lambda = \prod^*_{x_\lambda} (H \cap
  G_\lambda^{x_\lambda}) \ast F_\lambda$ such that $x_\lambda\Theta=1$ for all $x_\lambda$, and for each $\lambda$ the
  $x_\lambda$ runs through a suitable set of representatives of double cosets $G_\lambda x
  H$ such that $G_\lambda H$ is represented by $1$. Furthermore, the $F_\lambda$ are free.
\end{theorem}

Obviously, the $H_\lambda$ in the conjecture satisfy the requirements of Higgins's
Theorem. Therefore to prove the conjecture, one ``only'' has to decompose the
$H_\lambda$ provided by Higgins's Theorem. Fortunately, this is easy to do for the
$H_\lambda$ in Higgins's proof: the Kuro\v s's Theorem just provides the right
decomposition. This is what we are going to do.

The relevant additional property of the $H_\lambda$ in Higgins's proof is that their
intersection with $G_\lambda$ is contained in a conjugate of $H_\lambda$, see
Lemma~\ref{th:generalized-Higgins}.  This is proved by some additional
arguments to Higgins's proof.

Therefore we recall briefly Higgins's proof of his theorems in
Section~\ref{sec:groupoid-method} as done in \cite{MR97f:20032}.  This will
make the proof of our main lemma understandable to the reader not familiar
with the groupoid proofs.

In Section~\ref{sec:proof-higg-conj} we formulate a stronger version of
Higgins's Theorem as our main lemma. Then we prove the lemma and Higgins's
conjecture.

\section{Groupoid method}
\label{sec:groupoid-method}

In this section we recall briefly Higgins's proof of Theorems~\ref{th:Kurosh}
and \ref{th:Grushko-Higgins}.  See \cite[Chapter~14]{MR48:6288} for full
details.  We follow the discussion in \cite{MR97f:20032}.

Recall that a groupoid is a category in which every morphism is invertible.
Any group $G$ can be regarded as a groupoid with one object such that the
automorphism group of the object is $G$. If $H$ is a subgroup of $G$ then the
\emph{standard covering} $\gamma\colon \widetilde{G} \to G$ is a functor defined as
follows.  First, we define the groupoid $\widetilde{G}$.  The objects are the
right cosets of $H$ in $G$. Morphisms of $\widetilde{G}$ are $(N,g)\colon N \to Ng$
where $N$ is an arbitrary right coset and $g \in G$. Composition is defined by
$(N,g) \circ (Ng,h) := (N,gh)$. For example, the automorphism group of the coset
$H$ in $\widetilde{G}$ is isomorphic to $H$ via the map $(H,h) \mapsto h$.  Finally,
$\gamma$ is given by the formula $(N,g)\gamma:=g$.

We will also think of a groupoid as an oriented graph where the vertices are
the objects and the edges are the morphisms. In this sense, we will speak of
connected groupoids, trees and so on.

We will use free product of groupoids: let $G_\lambda$ be groupoids whose objects
are contained in a set $S$.  The free product of $G_\lambda$ is the groupoid
generated by the $G_\lambda$ in which only the necessary relations hold.  The
objects of the free product are the objects of all the $G_\lambda$ and hence is
contained in $S$.  Free product is similar to coproduct but some objects are
identified; that is the role of $S$. For example, every groupoid is the free
product of its connected components. If $S$ is a one-element set, this notion
is exactly the free product of groups. In the following, $S$ is always the set
of objects of the standard covering of a subgroup $H$ of a group $G$.

If $C$ is a connected groupoid and $\tau$ is a spanning tree then $\tau$ generates a
\emph{wide tree subgroupoid} $T$ of $C$, i.e.\ a subgroupoid in which there is
exactly one morphism between any two objects.  Then $C$ is isomorphic to the
free product and the direct product of $H$ and $T$, where $H$ is the
automorphism group of an object.  The canonical projection $\rho_\tau=\rho\colon C \cong H \ast T =
H \times T \to H$ is given by $(N,h) \mapsto \beta_N^{-1} \cdot h \cdot \beta_{Nh}$.  Here $\beta_N$ denotes
the unique isomorphism $N \to H$ in $T$.

Let us suppose now that $G=\prod^*_{\lambda \in \Lambda}G_\lambda$ is a free product of groups.  Now
the idea of the groupoid proofs of Kuro\v s's Theorem and Higgins's Theorem is that
the free decomposition of $G$ lifts to a free decomposition
\begin{equation}\label{eq:2}
  \widetilde{G} = \prod^*_{\lambda \in \Lambda} \widetilde{G_\lambda}
\end{equation}
where $\widetilde{G_\lambda}:=G_\lambda\gamma^{-1}$.  Using a suitable tree $\tau$, the projection $\rho$
maps this decomposition to a free decomposition of $H$ which will satisfy the
theorems.

In case of Higgins's Theorem (Theorem~\ref{th:Grushko-Higgins}), let $\Theta\colon
G=\prod^*_{\lambda \in \Lambda}G_\lambda \to B=\prod^*_{\lambda \in \Lambda}B_\lambda$ be a homomorphism such that $G_\lambda\Theta=B_\lambda$
and $H\Theta=B$.  We choose the tree $\tau$ such that the wide tree subgroupoid $T$
generated by $\tau$ is contained in $\ker \gamma\Theta$ (the full subgroupoid of
$\widetilde{G}$ consisting of morphisms mapped to identity by $\gamma\Theta$).  If $\tau$
is chosen with care, we will have
\begin{equation}
  \label{eq:1}
  H = \prod^* \underbrace{\widetilde{G_\lambda}\rho}_{H_\lambda}. 
\end{equation}
The condition $T \subseteq \ker \gamma\Theta$ will guarantee $H_\lambda\Theta \subseteq B_\lambda$.  See \cite{MR34:4344}
or \cite{MR48:6288} for more details.

In case of Kuro\v s's Theorem, we first decompose each $\widetilde{G_\lambda}$ into its
connected components $\widetilde{G_{\lambda,\mu}}$, which we further decompose to the
group $K_{\lambda,\mu}$ of one of its objects and a wide tree subgroupoid generated by
a tree $\tau_{\lambda,\mu}$. This leads to the free decomposition:
\begin{equation}\label{eq:3}
  \widetilde{G} = \prod^*_{\lambda,\mu} K_{\lambda,\mu} \ast F(\bigcup \tau_{\lambda,\mu})
\end{equation}
where $F(X)$ denotes the groupoid freely generated by the morphisms in $X$.
Note that a tree always generates a wide tree subgroupoid freely.

It is easy to see that $\bigcup \tau_{\lambda,\mu}$ is connected and hence contains a spanning
tree $\tau$. Now $\rho_\tau$ gives the free decomposition:
\begin{equation}
  \label{eq:4}
  H =  \prod^*_{\lambda,\mu} K_{\lambda,\mu}\rho_\tau \ast F(\bigcup \tau_{\lambda,\mu}\setminus\tau).
\end{equation}
It is easily seen that $K_{\lambda,\mu}\rho_\tau = H \cap G_\lambda^{x_{\lambda,\mu}}$ for some $x_{\lambda,\mu}$ and
$F(\bigcup \tau_{\lambda,\mu}\setminus\tau)$ is a free group so this gives the Kuro\v s decomposition of
$H$. Actually, $(H,x_{\lambda,\mu}^{-1})$ is the unique isomorphism in $F(\tau)$ between
$H$ and the object at which $K_{\lambda,\mu}$ is located. An easy argument, which we
omit, yields that the $x_{\lambda,\mu}$ form a set of  representatives of double
cosets $G_\lambda x H$.
If we have chosen $K_{\lambda,\mu}$ at the object $H$ whenever $\widetilde{G_{\lambda,\mu}}$
contains $H$, then the coset $G_\lambda H$ will be represented by $1$.

Higgins conjectured that both theorems can be proved using a common $\tau$, which
would lead to a common generalization of both theorems and their proofs.
In~\cite{MR97f:20032} it is shown that in general there is no tree $\tau$ which
is contained in both $\bigcup \tau_{\lambda,\mu}$ and $\ker \gamma\Theta$, so such a generalization
requires significant changes to the above proofs.

\section{Proof of Higgins's conjecture}
\label{sec:proof-higg-conj}

First we prove that the $H_\lambda$ in Theorem~\ref{th:Grushko-Higgins} has some
nice properties.
\begin{lemma}[Generalization of Higgins's Theorem]
  \label{th:generalized-Higgins}
  Suppose that a group homomorphism $\Theta\colon G=\prod^*_{\lambda \in \Lambda} G_\lambda \to B=\prod^*_{\lambda \in \Lambda} B_\lambda$
  between free products satisfies $G_\lambda\Theta=B_\lambda$ for all $\lambda \in \Lambda$. Let $H \subseteq G$ be a
  subgroup such that $H\Theta=B$.  Then $H = \prod^*_{\lambda \in \Lambda} H_\lambda$ such that for each
  $\lambda$ we have $H_\lambda\Theta=B_\lambda$, and there are representatives $\beta_{\lambda,\mu}$ of double
  cosets $G_\lambda x H$ such that $H \cap G_\lambda^{\beta_{\lambda,\mu}} \subseteq H_\lambda$ and $\beta_{\lambda,\mu} \in \ker \Theta$.
  Moreover, $G_\lambda H$ can be represented by $1$ for all $\lambda$ simultaneously.
\end{lemma}
\begin{proof}
We combine the ideas of Higgins's Theorem and Kuro\v s's Theorem from
Section~\ref{sec:groupoid-method} together.

We start with the proof of Higgins's Theorem and thus obtain a free
decomposition of $H$ into the $H_\lambda=\widetilde{G_\lambda}\rho_\tau$.  Now we use the proof of
Kuro\v s's Theorem for the tree $\tau$. We decompose $\widetilde{G_\lambda}$ into its
connected components $\widetilde{G_{\lambda,\mu}}$ and from every
$\widetilde{G_{\lambda,\mu}}$ we select the automorphism group $K_{\lambda,\mu}$ of an object.
Thus $K_{\lambda,\mu}\rho_\tau \subseteq H_\lambda$.  It is not obvious whether we obtain a free
decomposition like~\eqref{eq:4} but we still have $K_{\lambda,\mu}\rho_\tau = H \cap
G_\lambda^{\beta_{\lambda,\mu}}$ for some representatives $\beta_{\lambda,\mu}$ of double cosets $G_\lambda xH$.
We also have $(H,\beta_{\lambda,\mu}^{-1}) \in F(\tau) \subseteq \ker \gamma\Theta$.  Hence $\beta_{\lambda,\mu} \in \ker \Theta$.
The coset $G_\lambda H$ is represented by $1$ if we choose $K_{\lambda,\mu}$ at the object
$H$ when $\widetilde{G_{\lambda,\mu}}$ contains the object $H$.

Thus the $\beta_{\lambda,\mu}$ satisfy the lemma.
\end{proof}

This lemma together with Theorem~\ref{th:Kurosh} is enough to prove Higgins's
conjecture without using groupoids.

\begin{proof}[Proof of Theorem~\ref{th:Higgins-conjecture}]
The proof consists of two steps: first we decompose $H$ into $H_\lambda$ using
Lemma~\ref{th:generalized-Higgins} and, secondly,  Kuro\v s's Theorem will give
the required decomposition of $H_\lambda$.

By Lemma~\ref{th:generalized-Higgins}, we have a decomposition $H=\prod^*_{\lambda \in \Lambda}
H_\lambda$ with $H_\lambda\Theta \subseteq B_\lambda$ such that $H \cap G_\lambda^{\beta_{\lambda,\mu}}$ is contained in $H_\lambda$ for
some representatives $\beta_{\lambda,\mu}$ of double cosets $G_\lambda x H$ and $\beta_{\lambda,\mu}\Theta=1$. We
do not claim that these representatives give a Kuro\v s type decomposition; we
will modify them.

Applying Theorem~\ref{th:Kurosh} to $H_\lambda$ we obtain a Kuro\v s decomposition:
\begin{equation}
  \label{eq:5}
  H_\lambda = \prod^*_{\varepsilon,\delta} (H_\lambda \cap G_\varepsilon^\delta) \ast F_\lambda,
\end{equation}
where $F_\lambda$ is free. We claim that this decomposition is exactly the
decomposition of $H_\lambda$ the theorem requires. First of all, $F_\lambda$ will be the
free component. Now we examine the other components.

For every pair $\varepsilon,\delta$ in~\eqref{eq:5} $\delta$ lies in a double coset
$G_\varepsilon \beta_{\varepsilon,\mu} H$ i.e.\ 
\begin{equation}
  \label{eq:6}
  \delta = g \beta_{\varepsilon,\mu}h \quad \text{for some $g \in G_\varepsilon$ and $h \in H$}.
\end{equation}
Then we have
\begin{equation}
  \label{eq:7}
  H \cap G_\varepsilon^\delta = H \cap G_\varepsilon^{\beta_{\varepsilon,\mu} h} = \left( H \cap G_\varepsilon^{\beta_{\varepsilon,\mu}} \right)^h \subseteq H_\varepsilon^h.
\end{equation}
Therefore
\begin{equation}
  \label{eq:8}
  H_\lambda \cap G_\varepsilon^\delta = H_\lambda \cap \left( H \cap G_\varepsilon^\delta \right) \subseteq H_\lambda \cap H_\varepsilon^h = 
  \begin{cases}
    H_\lambda & \text{if $\varepsilon=\lambda$ and $h \in H_\lambda$,}\\
    1 & \text{otherwise.}
  \end{cases}
\end{equation}
In other words, $H_\lambda \cap G_\varepsilon^\delta$ is trivial unless $\varepsilon=\lambda$ and $\delta$ comes from a
double coset $G_\lambda \beta_{\lambda,\mu} H_\lambda$, and in this case $H_\lambda \cap G_\varepsilon^\delta=H \cap G_\lambda^\delta$.  So if
we denote by $\beta'_{\lambda,\mu}$ the representative of $G_\lambda \beta_{\lambda,\mu} H_\lambda$ occurring
in~\eqref{eq:5} then the free decomposition of $H_\lambda$ reduces to, after
omitting the components which~\eqref{eq:8} shows trivial:
\begin{equation}
  \label{eq:9}
  H_\lambda = \prod^*_{\mu} (H \cap G_\lambda^{\beta'_{\lambda,\mu}}) \ast F_\lambda.
\end{equation}
For each $\lambda$ the elements $\beta'_{\lambda,\mu}$ obviously form a set of double coset representatives
and $\beta'_{\lambda,\mu}\Theta \in (G_\lambda \beta_{\lambda,\mu} H_\lambda)\Theta =B_\lambda$. Since $G_\lambda\Theta=B_\lambda$, there are
elements $g_{\lambda,\mu} \in G_\lambda$ such that $g_{\lambda,\mu}\Theta=\beta'_{\lambda,\mu}\Theta$. Setting $x_{\lambda,\mu} :=
g_{\lambda,\mu}^{-1} \cdot \beta'_{\lambda,\mu}$, we have $x_{\lambda,\mu} \in \ker \Theta$ and
$G_\lambda^{\beta'_{\lambda,\mu}}=G_\lambda^{x_{\lambda,\mu}}$, hence
\begin{equation}
 \label{eq:10}
  H_\lambda = \prod^*_{\mu} (H \cap G_\lambda^{x_{\lambda,\mu}}) \ast F_\lambda.
\end{equation}
Moreover, the elements $x_{\lambda,\mu} \in G_\lambda \beta_{\lambda,\mu} H$ form a set of representatives
of double cosets $G_\lambda x H$. For the unique $\mu$ with $\beta_{\lambda,\mu}=1$, we have
$\beta'_{\lambda,\mu}=1$, and we may choose $g_{\lambda,\mu}=1$. This implies $x_{\lambda,\mu}=1$ and
hence $1$ occurs in the double coset representatives.
\end{proof}


\end{document}